\documentclass[12pt]{amsart}
\usepackage{amssymb}
\usepackage{amsmath}
\usepackage{pictex}
\newcommand\alp{t}
\newcommand\AND{\quad\mbox{and}\quad}

\newcommand\bet{s}
\newcommand\C{\mathbb C}
\newcommand\Cy{\mathsf{Cyc}}
\newcommand\BC{\mathcal B}
\newcommand\cc{\mathsf{c}}
\newcommand\CC{\mathcal C}
\newcommand\de{\delta}
\newcommand\DL{\mathsf{DL}}

\newcommand\FC{\mathcal F}
\newcommand\HC{\mathcal H}
\newcommand\JC{\mathcal J}
\newcommand\la{\lambda}
\newcommand\om{\omega}
\newcommand\MC{\mathcal M}
\newcommand\mm{\mathsf{m}}
\newcommand\N{\mathbb N}

\newcommand\R{\mathbb R}
\newcommand\rr{\mathsf{r}}
\newcommand\SC{\mathcal S}
\newcommand\spec{\operatorname{\sf spec}}
\newcommand\T{\mathbb T}
\newcommand{\uno}{{\boldsymbol 1}}
\newcommand\wh{\widehat}
\newcommand\wt{\widetilde}
\newcommand\Z{\mathbb Z}

\numberwithin{equation}{section}

\newtheoremstyle{mythm}
  {9pt}
  {9pt}
  {\itshape}
  {0pt}
  {\bfseries}
  {}
  { }
  {\thmnumber{(#2)}\thmname{ #1}\thmnote{ #3}}

\newtheoremstyle{mydef}
  {9pt}
  {9pt}
  {\normalfont}
  {0pt}
  {\bfseries}
  {}
  { }
  {\thmnumber{(#2)}\thmname{ #1}\thmnote{ #3}}

\theoremstyle{mythm}
\newtheorem{thm}[equation]{Theorem.}
\newtheorem{pro}[equation]{Proposition.}
\newtheorem{lem}[equation]{Lemma.}
\newtheorem{cor}[equation]{Corollary.}

\theoremstyle{mydef}
\newtheorem{dfn}[equation]{Definition.}
\newtheorem{exa}[equation]{Example.}

\newtheorem{rmk}[equation]{Remark.}

\begin{document}
\title{\large The spectrum of the averaging operator\\
on a network (metric graph)}
\author{\bf Donald I. CARTWRIGHT and Wolfgang WOESS}
\address{\parbox{0.7\linewidth}{School of Mathematics and Statistics,\\
University of Sydney, N.S.W.~2006, Australia}}
\email{donaldc@maths.usyd.edu.au}
\address{\parbox{0.7\linewidth}{Institut f\"ur Mathematik C,
Technische Universit\"at Graz,\\
Steyrergasse 30, A-8010 Graz, Austria}}
\email{woess@TUGraz.at}
\date{September 15, 2005}
\thanks{The second author was supported by a visiting scholarship at
the University of Sydney in April 2004, and the first author was
supported by a visiting professorship at TU Graz in June 2005}

\begin{abstract}
A network 
is a countable, connected graph $X$ viewed
as a one-complex, where each edge $[x,y]=[y,x]$  ($x,y \in X^0$, the vertex set)
is a copy of the unit interval within the graph's one-skeleton $X^1$
and is assigned a positive conductance $\cc(xy)$.
A reference ``Lebesgue'' measure on $X^1$ is built up by using
Lebesgue measure with total mass $\cc(xy)$ on each edge $[x,y]$.
There are three natural operators on $X$: the transition operator
$P$ acting on functions on $X^0$
(the reversible Markov chain associated with $\cc$),
the averaging operator $A$ over spheres of radius $1$ on $X^1$, and the
Laplace operator $\Delta$ on $X^1$ (with Kirchhoff conditions weighted by $\cc$
at the vertices). The relation between the $\ell^2$-spectrum
of~$P$ and the $H^2$-spectrum of~$\Delta$ was described by {\sc Cattaneo}~\cite{Cat}. In
this paper we describe the relation between the $\ell^2$-spectrum of
$P$ and the $L^2$-spectrum of $A$.
\end{abstract}

\maketitle

\markboth{\sf D. I. Cartwright and W. Woess}{\sf The averaging operator
on a network}
\baselineskip 15pt

\section{Introduction}\label{intro}

Let $X$ be a countable, connected graph with symmetric neighbourhood
relation $\sim$ and without loops and multiple edges. We shall view it as a
one-complex, where each edge is a (homeomorphic) copy of the unit interval
and edges are glued together at common endpoints (vertices). We write $X^0$ for
the vertex set and $X^1$ for the one-skeleton of $X$. Every point
in $X^1$ is of the form $(xy,\alp)$, the point at distance
$\alp$ from $x$ on the non-oriented edge $[x,y]=[y,x]$, where 
$0 \le \alp \le 1$, and $x, y \in X^0$, $x\sim y$. 
Thus, $(xy,0)=x$ and $(xy,\alp) = (yx,1-\alp)$.
In this way, the discrete graph metric $d(\cdot,\cdot)$ on the vertex set
(minimal length = number of edges of a connecting path) has a natural
extension to $X^1$.

We equip each edge $[x,y]$ with a positive \emph{conductance} $\cc(xy)=\cc(yx)$.
On $X^0$, we consider the discrete measure $\mm^0$, where
$\mm^0(x) = \sum_{y: y \sim x} \cc(xy)$. Our \emph{basic assumption} is that
$\mm^0(x) < \infty$ for all $x \in X^0$. On $X^1$, we introduce the continuous
weighted ``Lebesgue'' measure $\mm^1$ which at the point $(xy,\alp)$ is given by
$\cc(xy)\cdot d\alp$, if
$0 < \alp < 1$ (the vertex set has $\mm^1$-measure $0$).
The pair $(X,\cc)$, together with these measures, is called
a \emph{network}, or -- in the recent literature -- also \emph{metric graph}
or \emph{quantum graph}.

Associated with a network, there are three natural operators.

The first is the \emph{transition operator} 
$P$ acting on functions $g: X^0 \to \C$ by
\begin{equation}\label{eq:Paction}
Pg(x) = \frac{1}{\mm^0(x)} \sum_{y: y \sim x} \cc(xy)\,g(y)\,.
\end{equation}
The second is the \emph{Laplace operator} $\Delta$. It can be defined
via Dirichlet form theory, or by considering the space of all continuous
functions $F: X^1 \to \C$ which are twice differentiable in the
interior of each edge and satisfy the Kirchhoff equations
$$
\sum_{y: y \sim x} \cc(xy) \,F'(xy,0+) = 0 \quad \text{for all}\; x \in X^0\,.
$$
We then have
$$
\Delta F(xy,\alp) = F''(xy,\alp)\,,
$$
the 2nd derivative with respect to
$\alp \in (0\,,\,1)$, and $\Delta$ has to be closed suitably.
See e.g. {\sc Cattaneo}~\cite{Cat}, {\sc Solomyak}~\cite{Sol} or
{\sc Eells and Fuglede}~\cite{EeFu} for precise details. (The paper~\cite{Cat}
seems to have escaped the attention of most people working on metric graphs.)

The third operator is the \emph{averaging operator} $A$ over balls of
radius $1$. It acts on locally integrable functions $F: X^1 \to \C$ by
\begin{equation}\label{eq:Aaction}
AF(xy,\alp) =
\frac{1}{\mm^0(x)}\sum_{u \sim x} \cc(xu) \int_{0}^{1-\alp} F(xu,\bet)\,d\bet
+ \frac{1}{\mm^0(y)}\sum_{v \sim y} \cc(yv) \int_{0}^{\alp} F(yv,\bet)\,d\bet \,.
\end{equation}
In the \emph{regular} case, i.e., when
$\mm^0(\cdot)$ is constant, this is just the $\mm^1$-average of $F$ over the
ball with radius $1$ centered at $(xy,\alp)$.

Each of the three operators gives rise to a Markov process. For $P$,
this is the \emph{random walk} (reversible Markov chain) with discrete time
and state space $X^0$ whose transition probabilities are
$p(x,y) = \cc(xy)/\mm^0(x)$, if $y \sim x$, and $p(x,y) = 0$, otherwise.

The Laplace operator $\Delta$ is the infinitesimal generator of
\emph{Brownian motion} on the network.

The stochastic interpretation of $A$ is more similar to that of $P$.
Namely, $A$ governs the \emph{random walk}  with discrete time
and state space $X^1$, where at any time $n$, if the current position is
$(xy,\alp)$, the next step goes to a random point in $X^1$ at distance
at most $1$. The random choice depends on $\alp$, the conductance of $[x,y]$
and the edges incident with $[x,y]$.

These stochastic aspects are not at the heart of the present paper. What we
are interested in here is relation between the spectra of the operators
$A$ and $P$. {\sc Cattaneo}~\cite{Cat} has given a complete description of the $H^2$-spectrum
of $\Delta$ in terms of the $\ell^2$-spectrum of $P$. Our plan is to describe
the $L^2$-spectrum of $A$ in terms of the $\ell^2$-spectrum of $P$.

This refers to the (complex) Hilbert spaces $L^2(X^1,\mm^1)$ and
$\ell^2(X^0,\mm^0)$.
The inner product of the latter is given by
$$
\langle g_1, g_2 \rangle = \sum_{x \in X^0} g_1(x)\overline{g_2(x)}\,\mm^0(x)
\,,
$$
and it is well known and easy to check that $P$ is self-adjoint with 
$\|P\| \le 1$ on this space.
Analogously, the inner product on $L^2(X^1,\mm^1)$ is
$$
\langle F_1,F_2\rangle=\frac{1}{2}\sum_{x\in X^0}\sum_{y\in X^0:y\sim x} \cc(xy)
\int_0^1 F_1(xy,\alp) \overline{F_2(xy,\alp)} \ d\alp\,.
$$
The factor~$\frac{1}{2}$ occurs because
corresponding to each edge $[x,y]$, we get two equal terms on
the right, namely $\int_0^1F_1(xy,\alp)\overline{F_2(xy,\alp)}\ d\alp$
and $\int_0^1F_1(yx,\alp)\overline{F_2(yx,\alp)}\ d\alp$.
Again, it is straightforward to check that $A$ is self-adjoint with
norm bounded by $1$ on $L^2(X^1,\mm^1)$.

There is a large body of
literature on the spectrum of transition (resp. adjacency and discrete
Laplace) operators on finite graphs, see e.g. the books by 
{\sc Biggs}~\cite{Big}, {\sc Cvetkovi\'c, Doob and Sachs}~\cite{CvDoSa} 
and {\sc Chung}~\cite{Chu}. Transition operators on infinite graphs
are also very well studied objects, see e.g. the books by {\sc Soardi} 
\cite{Soa} and {\sc Woess}~\cite{Wbook}.
A lot is known about the $\ell^2$-spectrum of transition operators on 
various classes of infinite graphs,
see e.g. {\sc Mohar and Woess}~\cite{MoWo} for a general survey (up to 1989),
and the many more recent papers, mostly embedded into the context of
Markovian convolution operators on groups, of which we quote here only a few:
{\sc de la Harpe, Robertson, and Vallette~\cite{HaRoVa1},~\cite{HaRoVa2},
Cartwright~\cite{Car}, Grigorchuk and Zuk~\cite{GriZuk1},~\cite{GriZuk2},
Bartholdi and Woess~\cite{BaWo}.}

On the other hand, not much work has been done regarding the spectra
of averaging operators on networks, whence it appears to be useful to have a
recipe for translating the spectrum of $P$ into the spectrum of $A$.
Our main result is the following.
\begin{thm}\label{thm:mainresult}
The spectrum of~$A$ is
$$
\spec(A)= 
\bigl\{ 0 \bigr\} 
\cup \left\{\frac{\sin \om}{\om} : \om  \in \R \setminus \{0\}\,,\;
\cos \om \in \spec(P) \right\} \cup \bigl\{ 1 : 1 \in \spec(P) \bigr\}
\,.
$$
\end{thm}

Here, by ``$\cup\{ 1 : 1 \in \spec(P) \}$'' we mean
that $1$ is included in $\spec(A)$ if and only if $1\in\spec(P)$.
This theorem has the following obvious consequence.
\begin{cor}\label{cor:specrad}
Let $\rho = \rho(P)$ denote the spectral radius of~$P$. Then the spectral radius
of~$A$ is
$$
\rho(A) = \begin{cases} 1\,,&\text{if}\quad  \rho=1\,,\\
          \sqrt{1-\rho^2}\big/\arccos(\rho)\,,&\text{if}\quad \rho < 1\,.
          \end{cases}
$$
\end{cor}
Let $\spec_p(P)$ denote the point spectrum of~$P$, i.e., the
set of $\ell^2(X^0,\mm^0)$-eigenvalues of~$P$.

\begin{thm}\label{thm:pointspectrum} We have
$$
\spec_p(A)\setminus\{0\} = \bigl\{ 1 : 1 \in \spec_p(P) \bigr\}\cup
\left\{\frac{\sin \om}{\om} : \om  \in \R \setminus \{0\}\,,\;
\cos \om \in \spec_p(P) \right\} \,.
$$
Moreover,
$0\in\spec_p(A)$ unless $\mm^0(X^0) = \infty$ and $X$ is a tree
with the property that after removal of any edge, at least one of the two
connected components is recurrent.
\end{thm}

For the precise meaning of this last condition, see 
Definition~\ref{def:transience} and Proposition~\ref{pro:tree}
below.

The structure of this paper is as follows. In \S 2, we set up the basic tools
for relating $P$ and $A$. In \S 3, we study the contribution to the kernel
of $A$ that comes from flows in the network. In \S 4, we prove the above two
main theorems, and we also specify for finite graphs  
how one can obtain an orthonormal basis of $L^2(X^1,\mm^1)$ consisting
of eigenvectors (-functions) of $A$. In \S 5, we exhibit several examples.

\section{Interpolation of functions on the vertex set}\label{sect:interpol}

For $g\in\ell^2(X^0,\mm^0)$ and
$u\in L^2[0,1]$, define a function $F_{g,u}$ on~$X^1$ by
\begin{equation}\label{eq:fgudefn}
F_{g,u}(xy,\alp)=u(1-\alp)g(x)+u(\alp)g(y).
\end{equation}
Recall that $(xy,\alp) = (yx,1-\alp)$; the definition of $F_{g,u}$ is compatible
with this parametrization.
It is easy to check that $F_{g,u}\in L^2(X^1,\mm^1)$. In fact, it is routine to
calculate, for $g_1,g_2\in\ell^2(X^0,\mm^0)$ and $u_1,u_2\in L^2[0,1]$, that
\begin{equation}\label{eq:fguinnerprod}
\langle F_{g_1,u_1},F_{g_2,u_2}\rangle
=\langle g_1,g_2 \rangle \langle u_1,u_2 \rangle +
       \langle Pg_1,g_2 \rangle \langle u_1,Su_2 \rangle,
\end{equation}
where $S$ is defined in \eqref{eq:SJdefns} below, and 
$\;\langle u_1,u_2\rangle=\int_0^1u_1(\alp)\overline{u_2(\alp)}\,d\alp\;$
is the standard inner product on $L^2[0,1]$, while the inner products on
$L^2(X^1,\mm^1)$ and $\ell^2(X^0,\mm^0)$ are those defined in the introduction.

\begin{lem}\label{lem:action} The action of~$A$ on a function $F_{g,u}$ is
given by
$$
AF_{g,u}=F_{g,JSu}+F_{Pg,Ju}\,,
$$
where the operators $S$ and $J$ are given by
\begin{equation}\label{eq:SJdefns}
Su(\alp)=u(1-\alp) \AND Ju(\alp)=\int_0^\alp u(\bet)\;d\bet\,.
\end{equation}
\end{lem}

\begin{dfn}\label{M-def} Denote by~$\MC_0$ the linear span of the
functions $F_{g,u}$\,, where $g\in\ell^2(X^0,\mm^0)$, and where $u\in L^2[0,1]$.
Let $\MC$ denote the closure of~$\MC_0$ in~$L^2(X^1,\mm^1)$.
\end{dfn}

Lemma~\ref{lem:action} shows that $\MC_0$ is invariant under~$A$, and therefore
$\MC$ is too.

\begin{lem}\label{lem:orthogcompcond}
The orthogonal complement of~$\MC$ in~$L^2(X^1,\mm^1)$ consists
of the (equivalence classes of) square integrable functions $F:X^1\to\C$
such that for each $x\in X^0$,
\begin{equation}\label{eq:orthogcompcond}
\sum_{y: y\sim x}\cc(xy)\,F(xy,\alp)=0 \quad
\text{for almost all }\quad \alp\in[0,1].
\end{equation}
\end{lem}

\begin{proof} Suppose that $F\in\MC^\bot$. Then in particular,
$\langle F_{g,u},F\rangle=0$ for each $u\in L^2[0,1]$ and
for $g=\delta_x$, and for this~$g$, one calculates that
\begin{displaymath}
\langle F_{g,u},F \rangle =
\int_0^1 u(1-\alp)\overline{\sum_{y: y\sim x}\cc(xy)\,F(xy,\alp)}\,d\alp
\end{displaymath}
Since $u\in L^2[0,1]$ is arbitrary, \eqref{eq:orthogcompcond}~holds.

Conversely if \eqref{eq:orthogcompcond}~holds, then
$\langle F_{g,u},F\rangle=0$ for all $g\in\ell^2(X^0,\mm^0)$ of the form
$g=\delta_x$. By linearity, $\langle F_{g,u},F \rangle = 0$ if $g$ is finitely
supported, and using finite approximations, this implies that
$\langle F_{g,u},F \rangle = 0$
for all $g\in\ell^2(X^0,\mm^0)$ and $u\in L^2[0,1]$.
\end{proof}

\begin{cor}\label{cor:mperp}
The operator~$A$ leaves~$\MC$ and~$\MC^\bot$ invariant,
and is identically zero on~$\MC^\bot$.
\end{cor}

\begin{proof} If $F\in\MC^\bot$, then $AF=0$, as is immediate
from the definition of $A$ and Lemma~\ref{lem:orthogcompcond}.
\end{proof}

In the next Section~3, we shall give a complete description of
$\MC^\bot$ in terms of \emph{flows} in the
network, characterising, in particular, those networks  
for which $\MC^\bot = \{ 0 \}$. 

Let us record some elementary properties of the operators $J$ and ~$S$
arising in~\eqref{eq:SJdefns}.

\begin{lem}\label{lem:Sproperties}
The operator~$S$ satisfies $S^*=S$ and $S^2=I$. If $v,w\in L^2[0,1]$,
with $Sv=v$ (in which case we say that $v$ is \emph{even}) and $Sw=-w$
(in which case $w$ is called \emph{odd}), then $\langle v,w\rangle=0$. For any
$u\in L^2[0,1]$, we can write
\begin{equation}\label{eq:decomp}
u=v+w, \quad \text{where} \quad v=\frac{u+Su}{2} \AND w=\frac{u-Su}{2}\,,
\end{equation}
and then $Sv=v$, $Sw=-w$ (and so $\langle v,w\rangle=0$).
\end{lem}

\begin{lem}\label{lem:Jcompact} 
The operator $J$ is compact (in fact, Hilbert-Schmidt), but not normal.
Moreover, $J^*=SJS$, so that $JS$ and $SJ$ are self-adjoint
operators on~$L^2[0,1]$. For any $u\in L^2[0,1]$ we have
\begin{equation}\label{eq:sjplusjs}
JSu+SJu=\langle u,\uno\rangle\,\uno\,,
\end{equation}
where $\uno$ is the function taking the constant value~$1$ on~$[0,1]$.
\end{lem}

\begin{proof}For the compactness of~$J$, see {\sc Meise and Vogt}
\cite[Proposition~16.12 and Lemma~16.7(1)]{meisevogt} or {\sc Pedersen}
\cite[Proposition~3.4.16 and Lemma~3.4.5]{pedersen}. The other
assertions are easily checked.
\end{proof}

\begin{lem} Let $-1\le \la \le 1$. For $u,v\in L^2[0,1]$, write
\begin{displaymath}
\langle u,v\rangle_{\la}=\langle u,v\rangle+\la\langle u,Sv\rangle.
\end{displaymath}
Then
\begin{equation}\label{eq:innerprodineq}
(1-|\la|)\langle u,u \rangle \le \langle u,u \rangle_{\la}
\le (1+|\la|) \langle u,u \rangle\,.
\end{equation}
If $-1 < \la < 1$ then $\langle \cdot,\cdot \rangle_{\la}$ is an inner
product on~$L^2[0,1]$.
\\[4pt]
In the degenerate cases $\la = \pm 1$, we have
$\langle u,u \rangle_1 = 0 \iff Su = -u$, and
$\langle u,u \rangle_{-1} = 0 \iff Su = u$, respectively.
\end{lem}

\begin{proof}This is routine, using
$\langle u,Sv\rangle=\langle Su,v\rangle$.
\end{proof}

For $-1 < \la < 1$, we shall denote by $L^2_{\la}$ the space~$L^2[0,1]$
endowed with the inner product~$\langle \cdot,\cdot\rangle_{\la}$.
By~(\ref{eq:innerprodineq}), it is a Hilbert space.

\begin{lem}\label{lem:Jselfadjoint} Let $-1 < \la < 1$. Then
the operator $J_{\la} = JS+\la\, J$ is compact and self-adjoint
on the Hilbert space $L^2_{\la}$.
\end{lem}

\begin{proof} If $u,v\in L^2[0,1]$, then using $J^*=SJS$ and $(JS)^*=JS$,
\begin{displaymath}
\begin{aligned}
\langle J_{\la} u,v\rangle_{\la}
&= \bigl\langle \bigl(JS+ \la\, J\bigr)u,v\bigr\rangle +
\la\,\bigl\langle \bigl(JS + \la\, J\bigr)u,Sv\bigr\rangle\\
&=\langle JSu,v\rangle + \la\, \bigl(\langle Ju,v\rangle
+\langle JSu,Sv\rangle\bigr)+ \la^2\langle Ju,Sv\rangle\\
&=\langle u,JSv\rangle + \la\,
\bigl(\langle u,SJSv\rangle + \langle u,Jv\rangle\bigr)
+ \la^2\langle u,SJv\rangle\\
&=\langle u,J_{\la} v\rangle_{\la}.
\end{aligned}
\end{displaymath}
The compactness of~$J_{\la}$ follows from the compactness
of~$J$ on~$L^2[0,1]$, plus the fact that the norms of~$L^2[0,1]$
and~$L^2_{\la}$ are equivalent.
\end{proof}

It follows from~\cite[Proposition~16.2]{meisevogt} or
\cite[Theorem~3.3.8]{pedersen} that $L^2_\lambda$ has
an orthonormal basis consisting of eigenfunctions for~$J_{\la}$.
More explicitly:
\begin{lem}\label{lem:eigenfns} Let $-1 < \la < 1$ and set
$\om = \arccos \la \in (0\,,\,\pi)$. Then the functions
\begin{displaymath}
u_{\la,n}(\alp)=\frac{\sqrt{2}}{\sin \om}\, \sin\bigl((\om+2\pi n)\,\alp\bigr)
\,, \quad n\in\Z,
\end{displaymath}
form a complete orthonormal basis of~$L^2_{\la}$ consisting
of eigenvectors of~$J_{\la}$. In fact, $u_{\la,n}$ is
an eigenfunction for the eigenvalue
\begin{equation}\label{eq:mulambdan}
\mu_{\la,n} = \frac{\sin \om}{\om + 2\pi n}\,.
\end{equation}
\end{lem}
\begin{proof} Setting $u(\alp) = \sin(\vartheta \alp)$, 
where $\vartheta\ne0$, we compute
$$
\begin{aligned}
J_{\la} u(\alp) &=
\int_0^\alp\sin(\vartheta(1-\bet))\,d\bet +
        \la \int_0^\alp\sin(\vartheta\bet)\,d\bet\\
&=\frac{\sin \vartheta}{\vartheta}\,\sin(\vartheta\alp)
+\frac{\la - \cos\vartheta}{\vartheta}\,\bigl(1-\cos(\vartheta\alp)\bigr).
\end{aligned}
$$
So if $\vartheta$ is such that $\cos\vartheta=\la$, then
$J_{\la} u=\mu\, u$, where $\mu=\sin\vartheta/\vartheta$.
Taking $\vartheta=\om + 2\pi n$, we
see that $u_{\la,n}$ is an eigenfunction for the eigenvalue~$\mu_{\la,n}$.
Since these eigenvalues are distinct for distinct~$n$'s, the $u_{\la,n}$'s
are orthogonal in~$L^2_{\la}$. It is routine to check that
they are in fact orthonormal.

Suppose that $u\in L^2_{\la}$ and that
$\langle u,u_{\la,n}\rangle_{\la}=0$ for
all~$n\in\Z$. We claim that $u=0$.  Taking $v(\alp)=\sin(\vartheta\alp)$,
where $\cos\vartheta=\cos \om$, we find
that $\bigl(v+(\cos\om)\, Sv\bigr)(\alp) =
\sin(\vartheta)\cos\bigl(\vartheta(1-\alp)\bigr)$.
So from $\langle u,u_{\la,n}+\la\, Su_{\la,n}\rangle=0$ for
all~$n$ we find that
\begin{equation}\label{eq:orthog}
\int_0^1u(1-\alp)\,\cos\bigl((\omega+2\pi n)\alp\bigr)\,d\alp=0
\end{equation}
for all~$n$. Adding and subtracting (\ref{eq:orthog}), and
(\ref{eq:orthog}) with $n$ replaced by~$-n$, we find that for all $n\in\Z$,
\begin{displaymath}
\begin{aligned}
\int_0^1\sin(2\pi n\alp)\,\sin(\omega\alp)\,u(1-\alp)\,d\alp&=0\quad\text{and}\\
\int_0^1\cos(2\pi n\alp)\,\cos(\omega\alp)\,u(1-\alp)\,d\alp&=0.
\end{aligned}
\end{displaymath}
The first of these conditions implies that
$v(\alp)=\sin(\omega\alp)\,u(1-\alp)$ satisfies $v(1-\alp)=v(\alp)$ for almost
all~$\alp$, and the second condition implies that
$w(\alp)=\cos(\omega\alp)\,u(1-\alp)$ satisfies $w(1-\alp)=-w(\alp)$ for almost
all~$\alp$. That is, for almost all~$\alp$,
\begin{displaymath}
\sin\bigl(\omega(1-\alp)\bigr)\,u(\alp)=\sin(\omega\alp)\,u(1-\alp)
\AND \cos(\omega(1-\alp))\,u(\alp)=-\cos(\omega\alp)\,u(1-\alp).
\end{displaymath}
Multiplying the first of these equations by $\cos(\omega\alp)$ and the second
by $\sin(\omega\alp)$ and adding, we find that $u(\alp)=0$ almost everywhere.
Hence the family $\{u_{\la,n}:n\in\Z\}$ is a complete orthonormal basis for
$L^2_{\la}$.
\end{proof}

\section{Flows, and the space $\MC^\perp$}\label{sect:perp}

We now study in detail the space $\MC^\perp$ defined in 
Lemma~\ref{lem:orthogcompcond}; recall the defining 
relation~\eqref{eq:orthogcompcond}. 

Given our graph $X$, we consider the edge set $E=E(X)$ to be the set of ordered
pairs $xy$, where $x,y \in X^0$ and $x \sim y$. We set $\rr(xy) = 1/\cc(xy)$,
the \emph{resistance} of the edge $xy$. Let $\ell^2(E,\rr)$ be the
Hilbert space of all functions $\Phi: E(X) \to \C$ for which
$\langle \Phi, \Phi \rangle < \infty$, where the inner product is
$$
\langle \Phi, \Psi \rangle = \frac12 \sum_{x,y: x \sim y} 
\Phi(xy)\overline{\Psi(xy)} \, \rr(xy)\,.
$$

\begin{dfn}\label{dfn:flow} A \emph{flow} on the network $(X,\cc)$
is a function $\Phi \in \ell^2(E,\rr)$ such that
\begin{equation}\label{eq:kirchhoff}
\sum_{y:y \sim x} \Phi(xy)=0 \qquad \text{for all}\; x \in X\,.
\end{equation}
The flow is called \emph{odd,} if $\Phi(xy)=-\Phi(yx)$, and it is
called \emph{even}, if $\Phi(xy)=\Phi(yx)$. 
\end{dfn} 

Our definition requires, in particular, that 
$\langle \Phi,\Phi \rangle < \infty$. The latter number is often called
the \emph{energy} -- or, more appropriately, the \emph{power} -- of the flow
$\Phi$. In the literature, the term \emph{flow} usually applies to what
we call an \emph{odd} flow here. In this case, one may imagine each edge
$[x,y]$ as a tube with unit length and cross-section $\cc(xy)$, the tubes
are connected at the vertices of $X$, and the network of tubes is filled with
liquid. Then $\Phi(xy)$ is the amount of liquid per time unit that flows from
$x$ to $y$, whence $-\Phi(xy)=\Phi(yx)$ flows in the reverse direction. 
The condition \eqref{eq:kirchhoff} is Kirchhoff's law: the amount of
liquid per time unit that enters at any vertex coincides with the amount that
exits. In particular, our flows have no source or sink -- they are ``passive
flows''. In the above definition, even flows do not have such a nice
physical interpretation. We shall write $\JC^e$ and
$\JC^o$ for the (closed and orthogonal) subspaces of $\ell^2(E,\rr)$ consisting
of all even and odd flows on the network $(X,\cc)$, respectively.

\begin{rmk}\label{rem:bipartite}
A graph is called \emph{bipartite}
if we can partition its vertex set $X^0$ in two classes $C_1, C_2$
such that every edge has one endpoint in $C_1$ and the other in $C_2$.
Equivalently, this means that $X$ has no odd cycles (as defined below).
On a bipartite network, there is an obvious one-to-one correspondence between
odd and even flows:
\begin{equation}\label{eq:evenodd}
\Phi \in \JC^o \leftrightarrow \wt\Phi \in \JC^e\,,\quad
\text{where}\quad \wt\Phi(xy) = (-1)^i \Phi(xy)\,,\;\text{if}\; x \in C_i\;
(1=1,2)\,.
\end{equation}
\end{rmk}

Returning to $L^2(X^1,\mm^1)$, a function $F$ in that space is called 
\emph{even} if $F(xy,1-t)=F(xy,t)$, and \emph{odd} if $F(xy,1-t)=-F(xy,t)$,
for all $t \in [0\,,\,1]$ and each $xy\in E(X)$. Each $F$ has an orthogonal 
decomposition as a sum of its even and odd part.

It is straightforward to verify the following lemma.

\begin{lem}\label{lem:flowMperp} {\rm (a)} If $F \in \MC^\perp$ is even
(respectively, odd), and $u \in L^2[0,1]$, then 
$$
\Phi(xy) = \Phi_{F,u}(xy) 
= \cc(xy) \int_0^1 F(xy,t)u(t)\,dt
$$
defines an even (respectively, odd) flow with 
$\langle \Phi,\Phi \rangle \le \langle F,F \rangle \langle u,u \rangle\,$. 
\\[4pt]
{\rm (b)}  If $\Phi$ is an even flow and $u \in L^2[0,1]$ is even
(respectively, if $\Phi$ is an odd flow and $u \in L^2[0,1]$ is odd),
then
$$
F(xy,t) = \Phi(xy)u(t)/\cc(xy) 
$$ 
defines an even (respectively, odd) function in $\MC^\perp$ with
$\langle F,F \rangle = \langle u,u \rangle \langle \Phi,\Phi \rangle\,$. 
\end{lem}

The simple proof is left to the reader. Regarding (a), note that
when one of $F \in \MC^\perp$ and $u \in L^2[0,1]$ is even and the
other is odd, then $\Phi_{F,u}\equiv 0$. Thus, we may restrict to even $u$
when $F$ is even and to odd $u$ when $F$ is odd.  We
set $d^e = \dim \JC^e$ and $d^o = \dim \JC^o$ ($\le \infty$). 
In view of Lemma~\ref{lem:flowMperp}, the following is now the consequence of 
basic Fourier analysis.

\begin{pro}\label{pro:flowMperp} Let $\{ \Phi_m^e : 0 \le m < d^e \}$ and 
$\{ \Phi_m^o : 0 \le m < d^o \}$ be orthonormal
bases of the spaces $\JC^e$ and $\JC^o$, respectively. Then an orthonormal
basis of the subspace $\MC^\perp$ of $L^2(X^1,\mm^1)$ defined by  
\eqref{eq:orthogcompcond} is given by the set of all functions
$$
\begin{aligned}
&G_{m,0}^e(xy,t) = \frac{\Phi_m^e(xy)}{\cc(xy)}\,,\quad
G_{m,n}^e(xy,t) = \frac{\sqrt{2}\,\Phi_m^e(xy)\cos(2\pi n t)}{\cc(xy)}\,, \AND\\ 
&G_{m,n}^o(xy,t) = \frac{\sqrt{2}\,\Phi_m^o(xy)\sin(2\pi n t)}{\cc(xy)}\,,
\end{aligned}
$$
where $n \in \N=\{1,2,\ldots\}$, and $0 \le m < d^e$ or $0 \le m < d^o$, 
respectively.
\end{pro}

A \emph{cycle} in $X$ is a sequence $c=[x_0,\ldots,x_n]$
($n\ge3$) of vertices such that $x_0,\ldots,x_{n-1}$ are distinct,
$x_ix_{i+1} \in E(X)$ for $i=0,\ldots,n-1$, and $x_n=x_0$. 
Associated with $c$, there is the natural flow $\Phi_c \in \JC^o$ defined
by
\begin{equation}\label{eq:cycleflow}
\begin{aligned} 
&\Phi_c(x_ix_{i+1}) = 1\,,\quad  \Phi_c(x_{i+1}x_i) = -1 \quad
(i=0,\ldots,n-1)\,,\\
&\Phi_c(xy)=0\,, \quad  \text{if $xy$ is not an edge on $c$.}
\end{aligned}
\end{equation} 
We remark that our cycle $c=[x_0,x_1\ldots,x_n]$ has an orientation,
and that $\Phi_{c^*} = - \Phi_c$ when  $c^*=[x_n,x_{n-1}\ldots,x_0]$.

We now want to characterise those networks for which
$\MC^\perp = \{ 0 \}$. For this purpose we recall the following.

\begin{dfn}\label{def:transience}
The network $(X,\cc)$ is called \emph{transient,} if
$\sum_{n \ge 0} \langle P^n \delta_x, \delta_y \rangle < \infty$ for some 
(equivalently, for all) $x,y \in X$. Otherwise, the network is called 
\emph{recurrent.}
\end{dfn}

For the significance of this probabilistic notion, see e.g.
\cite{Wbook} or~\cite{Soa}. 

\begin{pro}\label{pro:tree}
One has $\MC^\perp=\{0\}$ if and only if\/ {\rm (i)}~$X$ is a tree
and\/ {\rm (ii)}~after removal of any edge, at least one of the two
connected components is recurrent as a subnetwork.
\end{pro}

\begin{proof} Suppose that $X$ has a cycle $c$.  Then by 
Lemma~\ref{lem:flowMperp}(b) we can use the odd flow
$\Phi = \Phi_c$ to construct a non-zero function in $\MC^\perp$.

Thus, $X$ has to be a tree if $\MC=L^2(X^1,\mm^1)$. Now suppose that
$X$ is a tree.

It follows from the \emph{flow criterion} for transience of networks
that on the tree $X$ there is a non-zero odd flow with finite power
if and only if there is an edge that disconnects $X$ into two transient
subtrees, see e.g.~\cite{Soa}, Theorems 3.33 and 4.20.
Thus, what is left is to show that on a tree, $\MC^\bot \ne \{0\}$ if and
only if there is a non-zero $\Phi \in \JC^o$.

Suppose that $\Phi$ is such a flow. Then Lemma~\ref{lem:flowMperp}(b) 
shows how one can construct a non-zero, odd function in $\MC^\perp$.
Conversely, suppose that $F \in \MC^\perp$ is non-zero. If the odd
(respectively, even) part $F^o$ (respectively, $F^e$) of $F$ is
non-zero then there must be an odd (respectively, even) function $u
\in L^2[0,1]$ such that the odd flow $\Phi_{F^o,u}$ (respectively,
even flow $\Phi_{F^e,u}$) defined in Lemma~\ref{lem:flowMperp}(a) is
non-zero. By Remark~\ref{rem:bipartite}, when $F^e\ne 0$, the even flow 
$\Phi_{F^e,u}$  can be transformed into a non-zero odd flow, 
since every tree is bipartite.
\end{proof}

Our final goal in this section is to describe how one finds (orthonormal) bases of 
$\JC^e$ and $\JC^o$, when $X$ is a \emph{finite} graph, in which case the flow
spaces do not depend on the specific conductances assigned to the edges. 
A \emph{spanning tree} of the graph $X$ is a subtree
$Y$ of $X$ which contains all vertices of $X$. It defines a subnetwork
$(Y, \cc_Y)$ whose conductance function $\cc_Y$ is the restriction of $\cc$
to $Y$. Recall that $E(X)$ consists of ordered pairs of adjacent vertices, 
that is, we have associated with each unoriented edge $[x,y]$ two oppositely
oriented edges $xy$ and $yx$.
It will be convenient to choose for each unoriented edge of $X$ one of its
endpoints as the \emph{initial} and the other as the \emph{terminal}
point. We write $\vec E(X)$ for the resulting set of oriented edges,
so that $E(X)$ is the disjoint union of the sets of ordered pairs 
$\{ xy : xy \in \vec E(X) \}$ and  $\{ xy : yx \in \vec E(X) \}$.
Also, we set $\vec E(Y) = \vec E(X) \cap (Y \times Y)$.

Consider an edge $xy \in \vec E(X) \setminus \vec E(Y)$. Adding this edge
to the tree, the new graph has precisely one cycle $c_{xy}=[x_0,\dots,x_k]$
($k \ge 3$) which is oriented such that $x=x_i$ and $y=x_{i+1}$ for some $i$.
We define
\begin{equation}\label{eq:cycles}
\Cy(X:Y) = \{ c_{xy} : xy \in \vec E(X) \setminus \vec E(Y) \}\,.
\end{equation}

The following is well known. 

\begin{lem}\label{lem:oddcycledecomp} Let $Y$ be a spanning tree of the finite
graph $X$. Then the set of flows $\{ \Phi_c : c \in \Cy(X:Y) \}$ 
is a basis of $\JC^o$.
Every odd flow $\Phi$ in $X$ has the unique decomposition
$$
\Phi = \sum_{xy \in \vec E(X) \setminus \vec E(Y)} 
\Phi(xy) \cdot \Phi_{c_{xy}}\,.  
$$
\end{lem}

\begin{proof} The function 
$$
\Psi=\Phi - 
\sum_{xy \in \vec E(X) \setminus \vec E(Y)} \Phi(xy) \cdot \Phi_{c_{xy}}
$$ 
vanishes on all edges in $E(X) \setminus E(Y)$.
Thus, $\Psi$ defines an odd flow in the finite tree $Y$, whence $\Psi \equiv 0$.
Linear independence of the $\Phi_c$, $c \in \Cy(X:Y)$, is immediate.
\end{proof}

If $X$ is finite and bipartite, then all cycles are even (have even length),
and \eqref{eq:evenodd} implies that the set of even flows
$$
\{ \wt\Phi_c : c \in \Cy(X:Y) \}
$$
is a basis of $\JC^e$. In general, the situation is
slightly more complicated. We decompose
\begin{equation}\label{eq:eocycles}
\Cy(X:Y) = \Cy^e(X:Y) \cup \Cy^o(X:Y)\,,
\end{equation}
where $\Cy^e(X:Y)$ consists of all even and $\Cy^o(X:Y)$ consists of all
odd cycles in $\Cy(X:Y)$. If $\Cy^o(X:Y) \ne \emptyset$, then we can choose
an (oriented) edge $x_0y_0 \in \vec E(X) \setminus \vec E(Y)$ such that 
$c_0=c_{x_0y_0}$ is an odd cycle. 
Now let $xy  \in \vec E(X) \setminus \vec E(Y)$ be any other edge such
that $c=c_{xy}$ is odd. We can define an associated even flow $\wt\Phi_{c_0,c}$
with $\wt\Phi_{c_0,c}(xy)=1$ by distinguishing the following two cases: 
(i) if $c$ and $c_0$ intersect, we can define $\Phi_{c_0,c}$
by suitably alternating the values $1$ and $-1$ on all edges of $c_0 \cup c$
whose endpoints do not lie on both $c$ and $c_0\,$;
(ii) If $c$ and $c_0$ do not intersect, we can define
$\Phi_{c_0,c}$ by suitably alternating the values $1$ and $-1$ on 
all edges of $c_0 \cup c$, and by suitably alternating the values $2$ and
$-2$ on all the edges of the unique path in the tree $Y$ that connects $c$
with $c_0$. The simple details are best understood by drawing a few figures.
 
\begin{lem}\label{lem:evencycledecomp} Let $Y$ be a spanning tree of the finite
graph $X$. Then the set of flows 
$$
\{ \wt\Phi_c : c \in \Cy^e(X:Y) \} \cup 
\{ \Phi_{c_0,c} : c \in \Cy^o(X:Y)\,,\; c \ne c_0 \} 
$$
(where $c_0 \in \Cy^o(X:Y)$, if the latter set is non-empty)
is a basis of $\JC^e$.
Every even flow $\Phi$ in $X$ has the unique decomposition
$$
\Phi = \sum_{\substack{\scriptstyle xy \in \vec E(X) \setminus \vec E(Y)\\
\scriptstyle c_{xy} \,\, \text{even}}}
\Phi(xy) \cdot \wt\Phi_{c_{xy}} +
\sum_{\substack{\scriptstyle xy \in \vec E(X) \setminus \vec E(Y)\\
\scriptstyle c_{xy} \,\, \text{odd},\ xy \ne x_0y_0}}
\Phi(xy) \cdot \wt\Phi_{c_0,c_{xy}}
\,.
$$
\end{lem}

\begin{proof} If $\Psi$ is the difference between $\Phi$ and the
sum on the right hand side, then $\Psi$ is an even flow on the graph
obtained from the tree $Y$ by adding the edge $x_0y_0$ (if $X$ has odd 
cycles) or just on the tree $Y$ itself (if $X$ is bipartite).
If $xy$ is an edge of that graph such that $y$ is the only neighbour
of $x$, \eqref{eq:kirchhoff} implies that $\Psi(xy)=0$, so that
$\Psi$ also is an even flow on the graph that remains after deleting
$x$ and the edge $xy$. Thus, after repeatedly ``chopping off'' finitely many
edges where $\Psi=0$, we are left with an even flow on the odd cycle 
$c_{x_0y_0}$, which must vanish on each edge. Thus, $\Psi \equiv 0$.
Once more, linear independence of the proposed basis is immediate.
\end{proof}

If $X$ is finite, the last two lemmas provide a simple algorithm
for finding bases of the (finite dimensional) spaces $\JC^o$ and
$\JC^e$, which can orthonormalized by
the Gram-Schmidt method. Then Proposition~\ref{pro:flowMperp}
leads to an orthonormal basis of the space $\MC^\perp \subset \ker A$.

\section{The spectral measure, and proof of the main results}\label{general}

Recall the Spectral Theorem for a normal operator~$T$ on a Hilbert
space~$\HC$ (see~\cite[Chapter 18]{meisevogt}
or~\cite[\S\S\,4.4, 4.5]{pedersen}, for example).
Let $\BC(\HC)$ denote the $C^*$-algebra of bounded linear operators
on~$\HC$. Let $C^*(T)$
denote the closure in~$\BC(\HC)$ of the space of polynomials
in~$T$ and~$T^*$. Then there is an isometric $*$-isomorphism
$\Phi:f\mapsto f(T)$ from the $C^*$-algebra $\CC\bigl(\spec(T)\bigr)$ of
continuous  functions on~$\spec(T)$ onto~$C^*(T)$. This isomorphism
maps the function $f(\lambda)\equiv\lambda^n$ to~$T^n$
for $n=0,1,\ldots$. Now let $\FC_\infty(K)$ denote
the $C^*$-algebra of bounded Borel measurable functions on the compact
set~$K \subset \C$.
Then there is a $*$-homomorphism
$\Psi:\FC_\infty\bigl(\spec(T)\bigr)\to\BC(\HC)$,
also written $f\mapsto f(T)$, which extends~$\Phi$, and which
is continuous in the following sense: if $(f_n)$ is a uniformly bounded
sequence of measurable functions on~$\spec(T)$ converging pointwise to
a function~$f$ on~$\spec(T)$, then $\langle f_n(T)x,y\rangle
\to\langle f(T)x,y\rangle$ for each $x,y\in\HC$.

For each Borel set $B\subset\spec(T)$, denote by~$E(B)$ the
operator $\Psi(\uno_B)$, where $\uno_B$ is the indicator
function of~$B$. Each $E(B)$ is a self-adjoint projection,
and the map $B\mapsto E(B)$ is called the {\it spectral
measure\/} of~$T$. For each $x,y\in\HC$, the map
$\mu_{x,y}:B\mapsto\langle E(B)x,y\rangle$ is a regular Borel
complex measure, and
\begin{displaymath}
\langle f(T)x,y\rangle
=\int_{\spec(T)}f(\lambda)\;d\mu_{x,y}(\lambda)
\end{displaymath}
for each $f\in\FC_\infty\bigl(\spec(T)\bigr)$. It is convenient to also
write $\int_{\spec(T)}f(\lambda)\,d\langle E_\lambda x,y\rangle$
for the integral on the right, interpreting the latter as a
Lebesgue-Stieltjes integral with respect to the function
$\la \mapsto \langle E_\lambda x,y\rangle =
\bigl\langle E\bigl((-\infty\,,\,\lambda]\bigr) x,y\bigr\rangle$.

Now suppose that $T$ is self-adjoint, so that $\spec(T)\subset\R$.
If $f$ is a bounded Borel measurable function defined on a Borel
set of~$\R$ containing~$\spec(T)$, then $f(T)$ is by definition
$f_{|\spec(T)}(T)$. For example, in Lemma~\ref{lem:cont} below, we
apply the above spectral theory to the
operators $A$ and~$J_{\la} = JS+\la\, J$, which are self-adjoint,
and of norm at most~$1$, so that their spectra are contained
in~$[-1,1]$. So if $f\in\FC_\infty([-1,1])$ we can
form the operators $f(A)$ and $f(J_{\la})$, acting on $L^2(X^1,\mm^1)$
and~$L^2_{\la}$, respectively, whenever $|\la| < 1$.

Similarly, it is convenient to define $E(B)=E\bigl(B\cap\spec(T)\bigr)$ for
any Borel subset $B$ of~$\R$. With this notation, if $\lambda\in\R$, then
$\lambda\in\spec(T)$ if and only if the operator
$E\bigl((\lambda-\epsilon\,,\,\lambda+\epsilon)\bigr)$ is non-zero for each
$\epsilon>0$. Equivalently, $\lambda \notin\spec(T)$ if and only if
$f(A)=0$, the zero operator, for
every continuous function $f$ supported in
$(\lambda-\epsilon\,,\,\lambda+\epsilon)$. Also,
$\lambda$ is an eigenvalue of~$T$ if and only if
$E(\{\lambda\})\ne 0$ (see~\cite[Lemma~18.5(3) and
Proposition~18.14]{meisevogt} or~\cite[Proposition~4.5.10]{pedersen}).

In the sequel, $E$ will always be the spectral measure of the operator $P$.

The following Perron-Frobenius-type proposition concerning $P$ can be found in
the literature in a few places, mostly under the assumption
that $\mm^0(\cdot)$ is bounded away from~0 on $X^0$ 
(in which case it becomes easier).
Since it appears not to be as well known as it should be, we include
its proof, whose first part is extrapolated from 
{\sc Kersting}~\cite[Lemma 3.1]{Ker}.

\begin{pro}\label{pro:eigenvalpmone} The operator~$P$ has eigenfunctions
in $ \ell^2(X^0,\mm^0)$ for the eigenvalue~$1$ if and only if
$\mm^0(X^0)<\infty$. In this case, the $1$-eigenspace $\HC_1$
consists of the constant functions on~$X^0$.

Furthermore, it has eigenfunctions in $ \ell^2(X^0,\mm^0)$ for the
eigenvalue~$-1$ if and only if $\mm^0(X^0)<\infty$ and the graph $X$ is
bipartite.  In this case, the $-1$-eigenspace $\HC_{-1}$
is spanned by the single function $\uno_{C_1} - \uno_{C_2}$,  where $C_1$
and $C_2$ are the two bipartite classes.
\end{pro}

\begin{proof} If $\mm^0(X^0)<\infty$, then the constant functions are 
in~$\ell^2(X^0,\mm^0)$ and are eigenfunctions of~$P$ for the eigenvalue~1.
If $\mm^0(X^0)<\infty$ and $X$ is bipartite, then $\uno_{C_1}-\uno_{C_2}$ is
an eigenfunction of~$P$ for the eigenvalue~$-1$.

Conversely, suppose that $g\in\ell^2(X^0,\mm^0)$ is nonzero, and 
$Pg=\la_0\, g$, where $\la_0\in\{-1,1\}$. We first show that $(X,\cc)$ must
be recurrent (see Definition~\ref{def:transience}). 
For $\lambda\in\R$, the operator $E(\{\lambda\})$ is the 
orthogonal projection of~$\ell^2(X^0,\mm^0)$ onto the 
$\lambda$-eigenspace of~$P$ (non-trivial if and only if $\la \in \spec_p(P)$). 
Then for each fixed $g_1,g_2\in\ell^2(X^0,\mm^0)$,
\begin{displaymath}
\begin{aligned}
\langle P^ng_1,g_2\rangle
&=\int_{[-1,1]}\lambda^n\,d\langle E_\lambda g_1,g_2\rangle\\
&=\langle E(\{1\})g_1,g_2\rangle
+(-1)^n\langle E(\{-1\})g_1,g_2\rangle+o(1)
\quad\text{as}\ n\to\infty
\end{aligned}
\end{displaymath}
by the Bounded Convergence Theorem. By hypothesis, $E(\{\lambda_0\}) g=g\ne0$, 
and there must be an~$x\in X$ so that 
$E(\{\lambda_0\})\delta_x\ne0$. So
\begin{displaymath}
\begin{aligned}
\langle P^{2n}\delta_x,\delta_x\rangle
&=\langle E(\{1\})\delta_x,\delta_x\rangle
+\langle E(\{-1\})\delta_x,\delta_x\rangle+o(1)\\
&=\|E(\{1\})\delta_x\|^2+\|E(\{-1\})\delta_x\|^2+o(1)
\end{aligned}
\end{displaymath}
tends to a nonzero limit as $n\to\infty$. 
Hence $\sum_{k=0}^\infty \langle P^{k}\delta_x,\delta_x\rangle=\infty$, 
and $(X,\cc)$ is recurrent.

Now $Pg=\la_0\, g$, and so $|g|=|Pg|\le P|g|$. Let $f=P|g|-|g|$. 
Since $\|\,|g|\,\|_2=\|g\|_2$,
\begin{displaymath}
\biggl\|\sum_{k=0}^{n-1}P^kf\biggr\|_2=\bigl\|P^n|g|-|g|\,\bigr\|_2
\le 2\|g\|_2\,,
\end{displaymath}
whence we have for each $x\in X^0$
\begin{displaymath}
\sum_{k=0}^{n-1}\langle P^{k}\delta_x,\delta_x\rangle f(x)
\le \biggl\langle\sum_{k=0}^{n-1}P^kf,\de_x \biggr\rangle
\le 2\sqrt{\mm^0(x)}\,\|g\|_2
\end{displaymath}
for each $x\in X^0$ and each integer $n \ge 1$. So if $f(x)>0$ for
some $x\in X^0$, then
\begin{displaymath}
\sum_{k=0}^\infty \langle P^{k}\delta_x,\delta_x\rangle 
\le 2\sqrt{\mm^0(x)}\,\|g\|_2 \big/ f(x)
<\infty\,,
\end{displaymath}
contradicting recurrence. So $f=0$. Therefore $|g|$ is a nonnegative
harmonic function, that is, $P|g|=|g|$, and so is constant by recurrence, 
see e.g. ~\cite[Theorem 1.16]{Wbook}.
Since the constant is nonzero, and $g\in\ell^2(X^0,\mm^0)$, we have 
$\mm^0(X^0)<\infty$. Now fix $x_0\in X^0$. Multiplying $g$ by a scalar,
we may assume that $g(x_0)=1$. Then $|g(y)|=1$ for all $y\in X$, and
\begin{displaymath}
\sum_{y: y \sim x_0} \frac{\cc(x_0,y)}{\mm^0(x_0)}|g(y)| =1 = g(x_0) 
= \la_0 \sum_{y: y \sim x_0} \frac{\cc(x_0,y)}{\mm^0(x_0)}g(y) 
= \biggl|\sum_{y: y \sim x_0} \frac{\cc(x_0,y)}{\mm^0(x_0)} \la_0 \,g(y)\biggr|.
\end{displaymath}
Hence equality holds in the triangle inequality, and therefore $\la_0\, g(y)=1$ 
for each $y\in X^0$ such that $y \sim x_0$. So if $\la_0=1$, the
connectedness of~$X$ implies that $g(y)=1$ for all~$y\in X^0$. 
If $\la_0=-1$, connectedness of~$X$ implies that 
$g(y)=(-1)^{\text{dist}(x_0,y)}$, and that $X^0$ is bipartite, 
with $C_1$ and~$C_2$ the sets of vertices at even
and at odd distance from~$x_0$, respectively.
\end{proof}

It follows from~\eqref{eq:fguinnerprod} that if $g_1,g_2\in\ell^2(X^0,\mm^0)$
are in two mutually orthogonal subspaces of~$\ell^2(X^0,\mm^0)$ which
are also $P$-invariant,
then $\langle F_{g_1,u_1},F_{g_2,u_2}\rangle=0$ for any
$u_1,u_2\in L^2[0,1]$. So if $\HC'$
denotes the orthogonal complement in~$\ell^2(X^0,\mm^0)$ of the sum of the
eigenspaces $\HC_1$ and~$\HC_{-1}$ (which are at most $1$-dimensional),
then the orthogonal decomposition
\begin{equation}\label{eq:hzerodecomp}
\ell^2(X^0,\mm^0)=\HC_1+\HC_{-1}+\HC'
\end{equation}
gives rise to a corresponding orthogonal decomposition of~$\MC$\,:
\begin{displaymath}
\MC=\MC_1+\MC_{-1}+\MC',
\end{displaymath}
where $\MC'$ is the closure of the
linear span of functions $F_{g,u}$, where $g\in\HC'$
and $u\in L^2[0,1]$, and $\MC_{\pm 1}$ are constructed analogously
from $\HC_{\pm 1}$.

\begin{lem}\label{lem:specialfns}
The subspaces $\MC_1$, $\MC_{-1}$ and~$\MC'$ are
invariant under~$A$. Let $u\in L^2[0,1]$. For $g\in\HC_1$,
$AF_{g,u}=F_{g,\langle u,\uno\rangle\uno}$.
For $g\in\HC_{-1}$, $AF_{g,u}=0$.
\end{lem}

\begin{proof} The invariance of the subspaces is immediate
from Lemma ~\ref{lem:action}. If $g\in\HC_1$, then
write $u=v+w$ as in~\eqref{eq:decomp}. Then $F_{g,w}=0$ because
$g$ is constant. Thus $F_{g,u}=F_{g,v}$, and $AF_{g,u}$ equals
\begin{displaymath}
AF_{g,v}=F_{g,JSv}+F_{Pg,Jv}=F_{g,JSv}+F_{g,Jv}=
F_{g,JSv+Jv}=F_{g,\langle v,\uno\rangle\uno}
=F_{g,\langle u,\uno\rangle\uno}
\end{displaymath}
because $JSv+Jv=2Jv$, and applying~\eqref{eq:decomp} to~$2Jv$
in place of~$u$, we see that the ``$v$-component" of~$2Jv$
is $Jv+SJv=JSv+SJv=\langle v,\uno\rangle\uno$ by
\eqref{eq:sjplusjs}.

Similarly, if $g\in\HC_{-1}$, then write $u=v+w$ as
in~\eqref{eq:decomp}. Then $F_{g,v}=0$
because $g(y)=-g(x)$ for each edge~$xy$.
Thus $F_{g,u}=F_{g,w}$, and $AF_{g,u}$ equals
\begin{displaymath}
AF_{g,w}=F_{g,JSw}+F_{Pg,Jw}=F_{g,JSw}+F_{-g,Jw}=
F_{g,JSw-Jw}=F_{g,-2Jw}=0
\end{displaymath}
because applying~\eqref{eq:decomp} to~$2Jw$
in place of~$u$, we see that the ``$w$-component" of~$2Jw$
is $Jw-SJw=-(JSw+SJw)=-\langle w,\uno\rangle\uno=0$ by
\eqref{eq:sjplusjs}.
\end{proof}

The following is one of our main tools for linking the spectra of $P$ and $A$.

\begin{pro}\label{pro:operatorpowerdecomp} Suppose that
$g_1,g_2\in\ell^2(X^0,\mm^0)$ and $u_1,u_2\in L^2[0,1]$. Then
for $n=0,1,\ldots$,
\begin{equation}\label{eq:operatorpowerdecomp}
\langle A^nF_{g_1,u_1},F_{g_2,u_2}\rangle
=\int_{\spec(P)}\langle J_{\la}^n u_1,u_2 \rangle_\lambda\;
d\langle E_\lambda g_1,g_2\rangle\,.
\end{equation}
\end{pro}

\begin{proof} The proof is by induction. By the Spectral Theorem for~$P$,
\begin{equation}\label{eq:pspectralthm}
\langle P^k g_1,g_2 \rangle =
\int_{\spec(P)} \lambda^k\, d\langle E_\lambda g_1,g_2 \rangle
\quad\text{for}\ k=0,1,\ldots.
\end{equation}
Case $n=0$ of~\eqref{eq:operatorpowerdecomp}
is immediate from~\eqref{eq:fguinnerprod}, together with the cases $k=0,1$
of~\eqref{eq:pspectralthm}.

Assume \eqref{eq:operatorpowerdecomp}~holds for $n$. By Lemma
\ref{lem:action}, $A^{n+1}F_{g_1,u_1}=A^n(F_{g_1,JSu_1}+F_{Pg_1,Ju_1})$.
This and the induction hypothesis yield
\begin{displaymath}
\begin{aligned}
\langle A^{n+1}F_{g_1,u_1},F_{g_2,u_2}\rangle
&=\langle A^nF_{g_1,JSu_1},F_{g_2,u_2}\rangle
+\langle A^nF_{Pg_1,Ju_1},F_{g_2,u_2}\rangle\\
&=\int_{\spec(P)} \Bigl(
                        \langle J_{\la}^n JSu_1 , u_2  \rangle_{\la}
                      + \langle J_{\la}^n Ju_1,u_2  \rangle_{\la}
                  \Bigr) \;d\langle E_\lambda g_1,g_2\rangle\\
&=\int_{\spec(P)}\langle J_{\la}^{n+1}u_1,u_2
                 \rangle_{\la}\;d\langle E_\lambda g_1,g_2\rangle
\end{aligned}
\end{displaymath}
because $\int_{\spec(P)}f(\lambda)\, d\langle E_\lambda Pg_1,g_2 \rangle
=\int_{\spec(P)} \lambda f(\lambda)\, d\langle E_\lambda g_1,g_2\rangle$
for any continuous function $f(\lambda)$ on~$[a,b]$, as
on sees by first taking $f(\lambda)=\lambda^k$, for $k=0,1,\ldots$,
and using~\eqref{eq:pspectralthm}. This completes the induction step.
\end{proof}

Let $u_1, u_2 \in L^2[0,1]$, and let $f: [-1\,,\,1] \to \C$ be bounded and
Borel measurable.
If $-1 < \la < 1$ then the function $\tilde f = \tilde f_{u_1,u_2}$,
\begin{equation}\label{eq:tildef}
\tilde f(\la) = \langle f(J_{\la})u_1, u_2 \rangle_{\la}
\end{equation}
is well defined, since $J_{\la}$ is selfadjoint on the Hilbert space
$L^2_{\la}$ by Lemma~\ref{lem:Jselfadjoint}.
When $\la = \pm 1$ then we cannot speak of $f(J_{\la})$, because
$\langle\cdot,\cdot\rangle_{\pm 1}$ is degenerate and we have no Hilbert space
on which $J_{\la}$ is selfadjoint.

We set $\spec(P)' = \spec(P) \cap (-1\,,\,1)\,$. Recall the definition
of the subspace $\HC'$ of $\ell^2(X^0,\mm^0)$ in~\eqref{eq:hzerodecomp}.

\begin{lem}\label{lem:cont}
If $u_1, u_2 \in L^2[0,1]$, and $f: [-1\,,\,1] \to \C$ is continuous,
then the function $\tilde f$ defined on $(-1\,,\,1)$
by~\eqref{eq:tildef}, is continuous, and
\begin{equation}\label{eq:contfn}
\langle f(A)F_{g_1,u_1},F_{g_2,u_2}\rangle
=\int_{\spec(P)'} \tilde f(\la)\;d\langle E_\lambda g_1,g_2\rangle
\end{equation}
for all $g_1\in\HC'$ and $g_2 \in \ell^2(X^0,\mm^0)$. If $f$ is the indicator
function $\uno_{\{\mu\}}$ of a singleton, then
$\tilde f$ is bounded and Borel measurable on~$(-1\,,\,1)$,
and \eqref{eq:contfn}~holds.
\end{lem}

\begin{proof}[Proof of Lemma~\ref{lem:cont}] First of all,
notice that when $g_1\in\HC'$, the integral in
Proposition~\ref{pro:operatorpowerdecomp} can be regarded as 
over~$\spec(P)'$. For $\HC_1$, $\HC_{-1}$
and $\HC'$  are the images of $E(\{1\})$, $E(\{-1\})$ and
$E\bigl((-1,1)\bigr)$, respectively. Hence
$E(B)g_1=E(B)E\bigl((-1,1)\bigr)g_1=E\bigl(B\cap(-1,1)\bigr)g_1=0$
for any $B\subset\{-1,1\}$. 

Next, for arbitrary $h: (-1\,,\,1) \to \C$ which is bounded and
Borel measurable, we have
\begin{equation}\label{eq:unihbd}
|\tilde h(\lambda)| \le
\|h(J_{\la})\|_\lambda \,\|u_1\|_\lambda\,\|u_2\|_\lambda
\le 2\|h\|_\infty\,\|u_1\|_2\,\|u_2\|_2\qquad \text{for all} \quad
\la \in (-1\,,\,1)\,,
\end{equation}
where $\|h(J_{\la})\|_\lambda$ is the operator norm on $L^2_{\la}$.

Now, given $f$ (continuous), choose a sequence $(p_n)$ of polynomials
such that $p_n \to f$ uniformly on $[-1\,,\,1]$. Clearly,
$\tilde p_n(\la)$ is a polynomial, and replacing
$h$ with $f-p_n$ in~\eqref{eq:unihbd}, we see that
$$
|\tilde f(\la) - \tilde p_n(\la)| \le
2\|f - p_n\|_\infty \, \|u_1\|_2 \, \|u_2\|_2
$$
for all $\la \in (-1\,,\,1)$. Consequently, $(p_n)$ 
converges uniformly to~$\tilde f$, and so the latter
is a continuous function on~$(-1\,,\,1)$.
By~\eqref{eq:operatorpowerdecomp}, we have
$$
\langle p_n(A)F_{g_1,u_1},F_{g_2,u_2}\rangle
=\int_{\spec(P)'} \tilde p_n(\la)\,
d\langle E_\lambda g_1,g_2\rangle\,,
$$
and letting $n\to \infty$, we see that this also holds for $f$ and $\tilde f$
in the place of $p_n$ and $\tilde p_n$, respectively.

Now let $f$ be the indicator function $\uno_{\{\mu\}}$ of a singleton.
Then there is clearly a uniformly bounded sequence $(f_n)$
of continuous functions converging pointwise to~$f$ on~$[-1,1]$.
Then for each fixed $\lambda\in(-1,1)$,
\begin{displaymath}
{\tilde f_n}(\lambda)=\langle f_n(J_\la)u_1,u_2\rangle_\lambda
\to\langle f(J_\la)u_1,u_2\rangle_\lambda
={\tilde f}(\lambda)
\end{displaymath}
by the Spectral Theorem applied to $J_\lambda$ on~$L^2_\lambda$. 
Hence $\tilde f$ is Borel measurable. Moreover, estimating as 
in~(\ref{eq:unihbd}), 
$|{\tilde f_n}(\lambda)|\le2\|f_n\|_\infty\|u_1\|_2\,\|u_2\|_2$,
and so the functions ${\tilde f_n}$ are uniformly bounded
on~$(-1,1)$. Hence (\ref{eq:contfn}) holds
for~$f$ by the Bounded Convergence Theorem, and by the Spectral
Theorem applied to~$A$, because it holds for each~$f_n$.
\end{proof}

\begin{rmk} The assertions in the lemma can be extended in two
directions, neither of which we need in the sequel:

(a) In Lemma \ref{lem:cont}, if $f$ is continuous, then a little 
more work in the proof shows that
$\tilde f$ has a continuous extension to~$[-1,1]$, and
\eqref{eq:contfn} is in fact valid for any $g_1\in\ell^2(X^0,\mm^0)$,
with $\spec(P)'$ replaced by~$\spec(P)$.

(b) Using~\cite[Lemma~18.1]{meisevogt}, it is easy to see that
the last statement in Lemma \ref{lem:cont} is valid for any
bounded measurable function~$f$.
\end{rmk}

\begin{proof}[Proof of Theorem~\ref{thm:mainresult}]
We must show that $\spec(A)$ equals~$\SC$, where
\begin{displaymath}
\SC=\{0\} \cup \{\mu_{\lambda,n}:\lambda\in\spec(P)',\; n\in\Z\} \cup
\{ 1 : 1 \in \spec(P) \}\,.
\end{displaymath}
First note that the set~$\SC$ is closed.
For suppose that $\mu\in\R$ is the limit of a sequence of points
$\mu_j=\mu_{\lambda_j,n_j}$ in~$\SC\setminus\{0,1\}$. If $|n_j|\to\infty$,
then $\mu_j\to0$, and so $\mu=0\in\SC$. So taking a
subsequence, we may assume that there is an $n\in\Z$ so that
$n_j=n$ for all~$j$. Since $\spec(P)$ is compact, taking a further
subsequence, we may suppose that $\lambda_j\to\lambda\in\spec(P)$.
If $\lambda\ne\pm1$, then $\mu=\lim_{j\to\infty}\mu_j=\mu_{\lambda,n}\in\SC$.
If $\lambda\in\{-1,1\}$, then $\mu_j\to0$ unless $\lambda=1$ and $n=0$,
in which case $\mu_j\to1\in\SC$.

To show that $\spec(A)\subset\SC$, suppose that 
$\mu_0\in [-1\,,\,1] \setminus\SC$. Then
$[\mu_0-\epsilon\,,\,\mu_0+\epsilon]\cap\SC=\emptyset$
for some $\epsilon>0$. If $\mu_0\ne1$, we may assume that 
$1\notin[\mu_0-\epsilon\,,\,\mu_0+\epsilon]$. Let $f$ be 
any continuous function
on $[-1\,,\,1]$ supported in $[\mu_0-\epsilon\,,\,\mu_0+\epsilon]$.
For any $\lambda\in\spec(P)'$, $f(J_{\la})=0$
because $f(\mu)=0$ for all $\mu\in\spec(J_{\la})$ ($\subset \SC$,
by Lemma~\ref{lem:eigenfns}). By Lemma~\ref{lem:cont},
$\langle f(A)F_{g_1,u_1},F_{g_2,u_2}\rangle=0$
for all $g_1\in\HC'$, $g_2\in\ell^2(X^0,\mm^0)$ and $u_1,u_2\in L^2[0,1]$.
Therefore $f(A)=0$ on the $A$-invariant subspace $\MC'$ of~$L^2(X^1,\mm^1)$.
By Corollary~\ref{cor:mperp} and Lemma~\ref{lem:specialfns}, $A=0$ 
on~$\MC^\perp$ and~$\MC_{-1}$, and is an orthogonal projection
on~$\MC_1$, and since $f(0)=0$, we get $f(A)=0$ on $\MC_{-1}+\MC^\perp$.
If $\mu_0\ne1$, then $f(A)=0$ on~$\MC_1$ because
then $f(0)=f(1)=0$, as arranged above. If $\mu_0=1$, then $1\notin\spec(P)$ 
and $\MC_1=\{0\}$. So in both cases, $f(A)=0$ on the whole 
of~$L^2(X^1,\mm^1)$ for every continuous
$f$ supported in $[\mu_0-\epsilon\,,\,\mu_0+\epsilon]$. Therefore
$\mu_0 \notin \spec(A)$, and $\spec(A)\subset\SC$.

\smallskip

For the reverse inclusion, suppose that $\mu_0\in\SC$.

Consider first  $\mu_0 \ne 0,1$. Then $\mu_0=\mu_{\lambda_0,n_0}$
for some $\la_0 \in \spec(P)'$ and $n_0 \in \Z$. Suppose that
$0<\epsilon<1-|\la_0|$. By~\cite[Proposition~4.4.5]{pedersen},
there is a $g\in\ell^2(X^0,\mm^0)$ so
that $\|Pg-\la_0g\|<\epsilon\|g\|$. Let $u=u_{\la_0,n_0}$. Then
\begin{displaymath}
AF_{g,u}-\mu_0F_{g,u}=F_{g,JSu}+F_{Pg,Ju}-\mu_0F_{g,u}=F_{Pg-\la_0g,Ju},
\end{displaymath}
and this has norm at most $2\|Pg-\la_0g\|\,\|u\|\le2\epsilon\|g\|\,\|u\|$,
by~\eqref{eq:fguinnerprod}. Again by~\eqref{eq:fguinnerprod},
\begin{displaymath}
\|F_{g,u}\|^2\ge(1-|\la_0|-\epsilon)\|g\|^2\|u\|^2
\end{displaymath}
because $|\langle Pg,g\rangle|\le(|\la_0|+\epsilon)\|g\|^2$. Therefore
\begin{displaymath}
\|AF_{g,u}-\mu_0F_{g,u}\|^2\le
\frac{4\epsilon^2}{1-|\la_0|-\epsilon}\|F_{g,u}\|^2,
\end{displaymath}
and since $F_{g,u}\ne0$, it follows that $\mu_0\in\spec(A)$.

Next, suppose that $\mu_0 = 1 \in \SC$. This can only happen when $1
\in \spec(P)$. If $1$ is not an isolated point of $\spec(P)$ then
there is a sequence $(\la_n)$ in $\spec(P)$ such that $\la_n \to 1$
from below. But then we just showed that $\mu_{\la_n,0} \in \spec(A)$,
and $\mu_{\la_n,0} \to 1$. Therefore, $1 \in \spec(A)$. If $1$ is an
isolated point of $\spec(P)$, then it must be an eigenvalue
(\cite[Proposition~18.14(3)]{meisevogt} or
\cite[Proposition~4.4.5]{pedersen}), and by
Proposition~\ref{pro:eigenvalpmone}, the constant function $g\equiv 1$
is an associated eigenfunction in $\ell^2(X^0,\mm^0)$. But then the
constant function $F \equiv 1$ on $X^1$ is in $L^2(X^1,\mm^1)$, whence
$1 \in \spec(A)$.

Finally, consider $\mu_0=0 \in \SC$. If $\spec(P)$ contains some
$\la \in (-1\,,\,1)$ then $\spec(A) \ni \mu_{\la,n} \to 0$ as
$|n| \to \infty$, whence $0 \in \spec(A)$. Otherwise,
$\spec(P) \subset \{-1\,,\,1\}$, so that by Proposition~\ref{pro:eigenvalpmone},
the space $\ell^2(X^0,\mm^0)=\MC_1+\MC_{-1}$ is at most 2-dimensional, 
which can happen only when $X$ has exactly two vertices and one edge. But in
this case, $\spec(A) = \{0, 1\}$, by Lemma~\ref{lem:specialfns}.
\end{proof}

\begin{proof}[Proof of Theorem~\ref{thm:pointspectrum}]
We first show that
\begin{equation}\label{eq:showthat}
\spec_p(A) \setminus \{ 0, 1 \}
= \bigl\{ \mu_{\la,n} : n \in \Z\,,\; \la \in \spec_p(P)\setminus\{-1,1\}
  \bigr\}.
\end{equation}
If $\lambda\in\spec_p(P)\setminus\{-1,1\}$ and $n\in\Z$, let
$g\in\ell^2(X^0,\mm^0)$ be  non-zero,
and satisfy $Pg=\lambda g$. Let $u=u_{\lambda,n}$. Then
$AF_{g,u}=\mu_{\lambda,n}F_{g,u}$ and $F_{g,u}\ne0$.
So $\mu_{\lambda,n}\in\spec_p(A)$.

Since $\om = \arccos \la \in (0\,,\,\pi)$, note that if $n\in\Z$ then,
\begin{displaymath}
0 < |\mu_{\lambda,n}| = \Bigl|\frac{\sin \om}{\om+2\pi n}\Bigr| \le
\frac{\sin\omega}{\omega}<1
\end{displaymath}
So $\mu_{\lambda,n}\in\spec_p(A)\setminus\{0,1\}$.

On the other hand, let $\mu\in\spec_p(A)\setminus\{0,1\}$.
Since $\mu\in\spec(A)$, we can write $\mu=\sin \om/\om$
for some $\om \in \R$ with $\la = \cos \om \in \spec(P)$, or equivalently,
$\mu=\mu_{\lambda,n}$ for some $n\in\Z$. Since $\mu\ne 0$,
a glance at the curve $\om \mapsto \sin \om/\om$ shows that the number
of solutions $\om$ to $\sin\om/\om = \mu$ is finite. That is,
the number of pairs $(\lambda,n)\in\spec(P)\times\Z$ for which
$\mu_{\lambda,n}=\mu$ is finite.
Thus, the set $F_{\mu}$ of all $\lambda\in\spec(P)$ such that
$\mu=\mu_{\lambda,n}$ for some~$n\in\Z$ is a finite set.
Also, $F_{\mu} \subset \spec(P)'$, since $\mu \notin \{0,1\}$.

Let $f=\uno_{\{\mu\}}$. By assumption, $\mu\in\spec_p(A)$. So
$f(A)$, being the orthogonal projection onto the
$\mu$-eigenspace of~$A$, is nonzero.
But $f(A)F=0$ for all $F \in \MC_1 \cup \MC_{-1}\cup\MC^\bot$, because
$\mu\ne0,1$. So $f(A)$ cannot vanish on $\MC'$. Thus there exist 
$g \in \HC'$ and $u \in L^2[0,1]$ so that 
$\langle f(A)F_{g,u},F_{g,u}\rangle > 0$.
Then by (\ref{eq:contfn}), we must have
$\langle f(J_\lambda)u,u\rangle_\lambda > 0$ for some
$\lambda\in\spec(P)'$.  For any such~$\lambda$, $f(J_\lambda)\ne0$,
and so $\mu=\mu_{\lambda,n}$
for some~$n$. Therefore $\lambda\in F_\mu$.
So the integrand in~(\ref{eq:contfn}) is nonzero only for
$\lambda\in F_\mu$, at most, and so the measure
$B\mapsto\langle E(B)g,g\rangle$ must assign nonzero measure to~$F_\mu$,
and therefore to~$\{\lambda\}$ for some $\lambda \in F_{\mu}$.
This $\lambda$ must be in $\spec_p(P)$, and $\mu=\mu_{\lambda,n}$
for some~$n\in\Z$. This completes the proof of~\eqref{eq:showthat}.

Now suppose that $1\in\spec_p(P)$. Then $0,1\in\spec_p(A)$ because
there are eigenfunctions of~$A$ in~$\MC_1$ for both~0 and~1,
by Lemma~\ref{lem:specialfns}. If $1\not\in\spec_p(P)$,
then $\HC'=\ell^2(X^0,\mm^0)$ and $\spec(P)'=\spec(P)$. So by
Lemma~\ref{lem:cont}, applied to $f=\uno_{\{0\}}$ and $f=\uno_{\{1\}}$,
we see that $f(A)F=0$ for any $F\in\MC$,
because 0 and~1 are not in the point spectrum of any $J_\lambda$,
$\lambda\in(-1,1)$. Of course $0\in\spec_p(A)$ if
$\MC\subsetneqq L^2(X^1,\mm^1)$, since $A$ is zero on~$\MC^\bot$.
This, together with Propositions~\ref{pro:tree} and~\ref{pro:eigenvalpmone},
proves the last statement of Theorem~\ref{thm:pointspectrum}.
\end{proof}

Theorem~\ref{thm:pointspectrum} and its proof applies, in particular,
to the case when the graph $X$ is finite. In this case, we also obtain
an orthonormal basis of the operator~$A$ by combining 
Theorem~\ref{thm:pointspectrum} with the following. In the proof 
of the next lemma, it is convenient to 
define $L^2_\la$ also when $\la \in \{-1,1\}$.
We define $L^2_1$ (respectively $L^2_{-1}$) to be the set of
$u\in L^2[0,1]$ such that $Su=u$ (respectively $Su=-u$). Note 
that for both these~$\lambda$'s, $\langle u,v\rangle_\la
=2\langle u,v\rangle$ for $u,v\in L^2_\la$, so that
$L^2_\la$ is a Hilbert space. Recall the 
definition~\ref{M-def} of $\MC_0$.

\begin{lem}\label{lem:finitecase} Suppose that $X$ is
finite, and let $g_1,\ldots,g_m$ be an orthonormal basis 
for $\ell^2(X^0,\mm^0)$
consisting of eigenfunctions of~$P$, with $Pg_j=\la_j\, g_j$
(with $-1 \le \la_j \le 1$) for each~$j$.
Then $\MC_0$ is closed, and consists of all functions
\begin{equation}\label{eq:functionsum}
F=F_{g_1,u_1}+\cdots+F_{g_m,u_m},
\end{equation}
where $u_j\in L^2[0,1]$ for each~$j$, and where
\begin{equation}\label{eq:conditions}
Su_j=u_j\;\;\text{if}\;\; \la_j=1,\quad
\text{and}\quad Su_j=-u_j\;\; \text{if}\;\; \la_j=-1.
\end{equation}
Moreover, the map $F\mapsto(u_1,\ldots,u_m)$ is a linear
isometry of $\MC_0$ onto the orthogonal direct sum
of the spaces $L^2_{\la_j}$, $j=1,\ldots,m$.
\end{lem}

\begin{proof} Each $g\in\ell^2(X^0,\mm^0)$ is a linear
combination of the $g_j$'s, and so each $F_{g,u}$ is a
sum~\eqref{eq:functionsum} of functions $F_{g_j,u_j}$.
If $\la_j=1$, then $g_j$ is constant on~$X^0$, and so
for each edge $xy$, and any $u\in L^2[0,1]$,
\begin{displaymath}
F_{g_j,u}(xy,\alp)=u(1-\alp)g_j(x)+u(\alp)g_j(y)
=(u(1-\alp)+u(\alp))g_j(x).
\end{displaymath}
Therefore, if $Su=-u$, we have $F_{g_j,u}\equiv 0$.
Applying~\eqref{eq:decomp} to $u=u_j$, we
see that $F_{g_j,u_j}=F_{g_j,v}+F_{g_j,w}=F_{g_j,v}$.
Similarly, if $\la_j=-1$, then
$F_{g_j,u_j}=F_{g_j,w}$, since $g_j(y)=-g_j(x)$ for any edge $xy$.
Therefore each $F_{g,u}$ can be written in the form~\eqref{eq:functionsum},
where \eqref{eq:conditions}~holds.

If $F$ is as in~\eqref{eq:functionsum}, then
$\langle F_{g_j,u_j},F_{g_k,u_k}\rangle=0$ if $j\ne k$, and so
\begin{displaymath}
\left\|\sum_{j=1}^mF_{g_j,u_j}\right\|^2
=\sum_{j=1}^m\Bigl(\langle g_j,g_j\rangle\langle u_j,u_j\rangle+
\la_j\, \langle g_j,g_j\rangle\langle u_j,Su_j\rangle\Bigr)
=\sum_{j=1}^m\langle u_j,u_j\rangle_{\la_j}.
\end{displaymath}
It follows that $\MC_0$ is isometric to the direct
sum of the Hilbert spaces $L^2_{\la_j}$,
recalling the special definition of~$L^2_\la$ made
above when $\la = \pm 1$. Therefore $\MC_0$ is
complete for its inner product, and so closed in~$L^2(X^1,\mm^1)$.
\end{proof}

\begin{cor}\label{cor:finitecase} Suppose that $X$ is
finite, and let $g_1,\ldots,g_N$ be an orthonormal basis for 
$\ell^2(X^0,\mm^0)$ consisting of eigenfunctions of~$P$, 
with $Pg_j=\la_j\, g_j$ (where $-1 \le \la_j \le 1$) for each~$j$.
We assume that $\la_1=1$ and $g_1$ is constant, and when $X$ is bipartite,
that $\la_N=-1$ and $g_N$ is constant on each bipartite class.
Then 
$$
\spec(A) = \spec_p(A) 
=\{ 0, 1 \} \cup \{ \mu_{\la_j,n} : |\la_j| < 1\,,\;n \in \Z \}\,,
$$
with $\mu_{\la_j,n}$ as in \eqref{eq:mulambdan}.
\\[5pt]
{\rm (i)} An orthonormal basis of the subspace $\MC$ of $L^2(X^1,\mm^1)$ consisting
of eigenvectors of $A$ is obtained as follows.
(Note that $\mm^1(X^1)=\mm^0(X^0)/2$.)\\[4pt]
{\rm (a)} For the eigenvalue $\mu=1$, the eigenspace is spanned by the function
$$
F^{(1,0)}(xy,t)=\frac{1}{\sqrt{\mm^1(X^1)}}\,.
$$
{\rm (b)} For the eigenvalue $\mu=0$, then 
the eigenspace is spanned by the functions
$$
F^{(1,n)}(xy,\alp)= \frac{\sqrt{2}}{\sqrt{\mm^1(X^1)}}\, \cos(2\pi n\alp)
\,,\quad n \in \N\,,
$$
if $X$ is not bipartite. If $X$ is bipartite, then the eigenspace is spanned
by the $F_{1,n}$ and the functions
$$
F^{(-1,n)}(xy,\alp)= (-1)^{i(x)}\frac{\sqrt{2}}{\sqrt{\mm^1(X^1)}}\, 
\sin(2\pi n\alp)\,,\quad n \in \N\,, 
$$
where $i(x)=1$ or $2$ according to whether $x$ lies in the bipartite class
$C_1$ or $C_2$.
\\[2pt]
{\rm (c)} For the eigenvalue $\mu$ with $0 < |\mu| < 1$, 
the eigenspace is spanned by the functions
$$
F^{(j,n)}(xy,\alp) = \frac{\sqrt{2}}{\sin\om_j}\,
          \Bigl(g_j(x)\sin\bigl((\om_j+2\pi n)(1-\alp)\bigr) + 
                g_j(y)\sin\bigl((\om_j+2\pi n)\alp\bigr) \Bigr)\,,
$$
where $(j,n) \in \{1,\ldots,N\}\times \Z$ is such
that $\mu_{\la_j,n}=\mu$ and $\om_j = \arccos \la_j$.
\\[5pt]
{\rm (ii)} An orthonormal basis of the subspace $\MC^\perp \subset \ker A$ 
is obtained via Proposition~\ref{pro:flowMperp} in combination with 
Lemmas~\ref{lem:oddcycledecomp} and~\ref{lem:evencycledecomp}.

\end{cor}

\section{Final remarks and examples}\label{final}

\begin{rmk}\label{remark:the}
One may ask why we call our operator $A$, defined in
\eqref{eq:Aaction}, \emph{the} averaging
operator, and not the one which takes the pure $\mm^1$-average over balls of
radius~$1$. The latter is given by
$$
\begin{aligned} \wt A F(xy,\alp) =
\ &\frac{1}{(1-\alp)\mm^0(x)+\alp\,\mm^0(y)}\\
& \times \left(
\sum_{u \sim x} c(xu) \int_{0}^{1-\alp} F(xu,\bet)\,d\bet
+ \sum_{v \sim y} c(yv) \int_{0}^{\alp} F(yv,\bet)\,d\bet\right) .
\end{aligned}
$$
The point is that unlike $A$, the latter operator does not enjoy a
nice and natural compatibility with the transition operator $P$ and
the Laplace operator on a network. Note, however, that $A = \wt A$
when $\mm^0(\cdot)$ is constant. This occurs, in particular, when
the graph $X$ is locally finite and regular, 
and $\cc(xy)=1$ for each edge $xy$, in which
case $P$ is called the \emph{simple random walk} (SRW) operator.
\end{rmk}

We now give a few examples of locally finite, regular graphs with 
conductances $\cc(xy)=1$ for each edge $[x,y]$, where the spectrum of 
$A$ can be computed via the (known) spectrum of~$P$.

\begin{exa}\label{example:Z} Equip the additive group $\Z$ of all integers with
the typical graph structure, where the edges are between $x$ and $x+1$,
$x \in \Z$. Then the SRW operator $P$, associated with
conductances~$\equiv 1$, is the convolution operator
$$
Pf = \varphi*f \quad \text{with} \quad
\varphi = \frac12 \bigl(\uno_{\{1\}} + \uno_{\{-1\}})\,.
$$
Computing the Fourier transform $\wh \varphi(\om) = \cos\om$,
$\om \in [0\,,\,2\pi]$, one finds the very well known fact that
$$
\spec(P) = [-1\,,\,1] \AND \spec_p(P)= \emptyset\,.
$$
On the other hand, the one-skeleton is the real line $\R$ (with the
integer points singled out as vertices), so that $A$ is the convolution
operator
$$
AF = \varphi*F\quad \text{with} \quad
\varphi = \frac12 \uno_{[-1\,,\,1]}\,.
$$
Since the Fourier transform of $\varphi$ is $\wh\varphi(\om) = (\sin \om)/\om$,
$\om \in \R$, one finds that
$$
\spec(A) = [(\sin \om^*)/\om^*\,,\,1] \AND  \spec_p(A) = \emptyset\,,
$$
where $\om^*$ is the smallest positive solution of the equation
$\tan \om = \om\,$, so that $(\sin \om^*)/\om^* =\cos \om^*\,$.
Numerical computation gives $\om^*=4.493409\dots$ and 
$\sin(\om^*)/\om^* = -0.217233\dots\,$.
\end{exa}

\begin{exa}\label{example:ZN} Equip the additive group $\Z_N =
\Z/(N\Z)=\{0,\ldots,N-1\}$ of integers modulo $N$ with the structure 
of a cycle, where $x \sim x+1$ (mod $N$) for $x \in \Z_N$. If each edge is assigned
conductance~$1$, then $\mm^0$ becomes twice the counting measure.
An orthonormal basis of eigenvectors of the SRW operator $P$ acting on 
$\ell^2(\Z_N,\mm^0)$
with associated eigenvalues is given by
$$
g_j(x) = (2N)^{-1/2}\,\exp(2\pi i\, x\, j/N) \AND \la_j = \cos(2\pi j/N)\,,\quad
j=0,\ldots, N-1\,,
$$
so that $\la_j = \la_{N-j}$ ($1 \le j \le N-1$) has multiplicity $2$
unless $j=N/2$ for even $N$.

The one-skeleton of $\Z_N$ can be identified with the torus, interpreted as
the compact additive group $\R/(N\R)$ of real numbers modulo $N$,
parametrized by the interval $[0\,,N)$, or equivalently, $(-N/2\,,\,N/2]$,
with Lebesgue measure.
Thus, as in Example~\ref{example:Z}, the averaging operator $A$ becomes the 
convolution operator $AF = \varphi*F$ with 
$\varphi = \frac12 \uno_{[-1\,,\,1]}\,$, but this time modulo $N$.
The Fourier transform (on $\Z$) of $\varphi$ is 
$\wh\varphi(n) = (\sin \varpi_n)/\varpi_n$, and 
$$
\spec(A) = \{(\sin \varpi_n)/\varpi_n : n \in \Z \}\,,\quad\text{where}\quad
\varpi_n = 2 \pi n/N\,.
$$
(We set $(\sin 0)/0=1$ by continuous extension.) 
The orthonormal basis of associated eigenfunctions is given by
$$
F_n(t) = N^{-1/2} \exp(i\, \varpi_n\, t)\,, \quad t \in [0\,,N)\,,\quad n \in \Z\,.
$$ 
Since the eigenvalues corresponding to $\varpi_n$ and $\varpi_{-n}$ coincide,
the eigenspace has dimension~2, unless $n$ is a multiple of $N$, or --
when $N$ is even -- a multiple of $N/2$. 

Relating this with 
Corollary~\ref{cor:finitecase}, we get the following: for $1 \le j < N/2$,
$\mu_{\la_j,n} =  \mu_{\la_{N-j},n} =(\sin \varpi_{j+nN})/\varpi_{j+nN}\,$,
and elementary computations yield for the functions $F^{(j,n)}$ and 
$F^{(N-j,n)}$ of  Corollary~\ref{cor:finitecase}(c) that
$F^{(j,n)}(xy,t) = F_{j+nN}(x+t)$ and $F^{(N-j,n)}(xy,t) = F_{-j-nN}(x+t)$, 
if $x \in \{0,\ldots,N-1\}$, $y=x+1$ (mod $N$), and $t \in [0\,,\,1]$.  

The eigenspace associated with the eigenvector~1 consists of course of the
constant functions.

Finally, $\ker A$ is spanned by the functions $F_{nN}(t)$, 
$n \in \Z\setminus\{0\}$, if $N$ is odd, and by the functions 
$F_{nN/2}(t)$, $n \in \Z\setminus\{0\}$, if $N$ is even.
Comparing with Corollary~\ref{cor:finitecase}, we have to consider 
the decomposition $\ker A = (\MC \cap \ker A) \oplus \MC^\perp$.
If $N$ is odd, then $\MC \cap \ker A$ is spanned by the functions
$F^{(1,n)}(xy,t)= \sqrt 2 \,\, \Re \bigl( F_{nN}(x+t)\bigr)$, and -- 
since $\JC^e = \{0\}$ and $\dim\JC^o=1$  -- the space $\MC^\perp$ is spanned 
by the functions $G^o_{0,n}(xy,t)= \sqrt 2 \,\, \Im\bigl( F_{nN}(x+t) \bigr)\,$,
$n \in \N\,$, where $\Re(\cdot)$ and $\Im(\cdot)$ denote real and imaginary
part. 

When $N$ is even, the situation is slightly more complicated. In this case,
$\JC^e$ is one-dimensional and spanned by the even flow with 
value $(-1)^x$ on the edge $[x,y]$ with $y=x+1$, $x=0, \ldots, N-1$.
The functions $F^{(-1,n)}$, $n \in \N$, of Corollary~\ref{cor:finitecase}(b)
are non-differentiable at the vertices, and the functions $G_{0,n}^e\,$, 
$n \in \N_0\,$, of Proposition~\ref{pro:flowMperp} are even discontinuous,
whence they have to be expressed as Fourier series in terms of 
the functions $F_{nN/2}$, $n \in \Z$.
\end{exa}

\begin{exa}\label{example:tree} Let $\T$ be the
homogeneous tree with degree $q+1$, where $q \ge 1$. It is well
known that the spectrum of the SRW operator on $\T$ is
$$
\spec(P) = [-\rho\,,\,\rho]\,,\quad \rho=\frac{2\sqrt q}{q+1}\,,
\AND \spec_p(P)= \emptyset\,,
$$
see~\cite[\S 7.B]{MoWo} and the references given there. 
Removal of any edge splits $\T$ into
two transient pieces. Therefore, in view of our theorems,
$\spec_p(A) = \{0\}\,$,
$$
\spec(A) = \left\{ \frac{\sin \om}{\om} :
  \om > 0\,,\;|\cos\om| \le \rho \right\}
  \cup \{ 0 \}\,,\AND
\rho(A) = \frac{q-1}{q+1}\bigg/\arccos\frac{2\sqrt q}{q+1}\,.
$$
The last formula for the spectral radius of $A$ was first found
by {\sc Saloff-Coste and Woess}~\cite{SalWoe} by a completely different
method, and indeed, the latter was the starting point for the present
investigation.

A closer look at $\spec(A)$ may be of interest. Let
$\om_{\pm \rho} = \arccos (\pm\rho)$, where $\rho = \rho(P)$, so that
$0 < \om_{\rho} < \pi/2$ and $\om_{-\rho} = \pi - \om_{\rho}$.
For $n \ge 0$, let $m_n$ and $M_n$ be the minimum and maximum, respectively,
of all numbers $|\sin\om|/\om$, where
$n\pi+\om_{\rho} \le \om \le n\pi+\om_{-\rho}$.
In particular, $M_0=\rho(A)$. Also, $M_n > M_{n+1}$ and $m_n > m_{n+1}$.
Then
$$
\spec(A) = \{0\} \,\cup\, \bigcup_{n\ge 0} [m_{2n}\,,\,M_{2n}]
\,\cup\, \bigcup_{n\ge 0} [-M_{2n+1}\,,\,-m_{2n+1}]\,.
$$
A closer analysis of these intervals shows that there is an $n_0$
such that $[m_{n}\,,\,M_{n}]$ and $[m_{n+2}\,,\,M_{n+2}]$ overlap
for all $n \ge n_0$, whence $\spec(A)$ is a union of finitely many
intervals, one of which contains $[-M_{n_0}\,,\,M_{n_0+1}]$ if $n_0$ is
odd, and $[-M_{n_0+1}\,,\,M_{n_0}]$ if $n_0$ is even.

When $q=2$, one can verify that $n_0=0$, whence
$\spec(A) = [-M_1\,,\,M_0]$, but as $q$ increases, $\spec(A)$ decreases
and becomes
a finite disjoint union of more and more intervals. Indeed, in the
limit, as $q \to \infty$, $\spec(A)$ tends to the set
$\{0\} \cup \{ (-1)^n/(\frac{2n+1}{2}\pi) : n \ge 0 \}$.

Finally, we determine $-M_1 = \min \,\spec(A)$. As in Example
\ref{example:Z}, let
$\om^* \in (\pi\,,\, \frac32\pi)$ be the smallest positive solution
of $\tan \om = \om$. Then
$$
-M_1 = \begin{cases} \dfrac{\sin \om^*}{\om^*} = \cos \om^*\,,
                     &\text{if} \; \pi + \om_{\rho} < \om^*\,, \\[10pt]
       \dfrac{\sin (\pi + \om_{\rho})}{\pi + \om_{\rho}}
    = -\dfrac{q-1}{q+1}\bigg/\left(\pi+\arccos\dfrac{2\sqrt q}{q+1}\right)\,,
        &\text{otherwise.}
       \end{cases}
$$
The first of these two cases holds precisely when
$\tan(\pi + \om_{\rho}) < \pi + \om_{\rho}$, or equivalently, when
$
\frac{q-1}{2\sqrt q} < \pi + \arccos\frac{2\sqrt q}{q+1}\,,
$
that is, for $q \le 82$ by numerical computation.

\smallskip

It seems unlikely that these results could have been found without
using the relation between the spectra of $P$ and $A$.
\end{exa}

\begin{exa}\label{example:DL} Let $\T_1$ and $\T_2$ be two
homogeneous trees with degrees $q+1$ and $r+1$, respectively.
In each of the trees, we choose a boundary point (end) and the associated
Busemann (horocycle) function $h: \T_i \to \Z$. The
\emph{Diestel Leader graph} $\DL(q,r)$ is the horocyclic product of
the two trees, i.e., the subgraph of their direct product,
$$
\DL(q,r) = \{ x_1x_2 \in \T_1 \times \T_2 : h(x_1)+h(x_2)=0 \}\,.
$$
See e.g.~\cite{BaWo} for a detailed description
and further references. In particular, in~\cite{BaWo} it is shown
that for the SRW operator $P$ on $\DL=\DL(q,r)$, the spectrum
$\spec(P) = [-\rho(P), \rho(P)]$ is pure point, i.e., it is the closure
of the point spectrum, and there is an orthonormal basis of
$\ell^2(\DL)$ consisting of finitely supported eigenfunctions of $P$.
One has
$$
\spec_p(P)= \left\{\rho(P) \cos \frac{m}{n}\pi : n \ge 2\,,\;1 \le m \le n-1
\right\}\,,\quad \rho(P)=\frac{2\sqrt{qr}}{q+r}\,.
$$
By our theorems, we can compute the spectrum of $A$, which is also pure point
and contains~0, since $\DL$ is not a tree.

\smallskip

In the specific case when $r=q$, $\DL(q,r)$ is a Cayley graph of the
\emph{lamplighter group} (wreath product) $\Z_q \wr \Z$, see again
\cite{BaWo} for details. In that case, the spectrum of $P$ had been
determined previously in~\cite{GriZuk1}, $\rho(P)=1$, and
the point spectrum of $A$ has the following particularly nice form:
$$
\spec_p(A) = \{0\}\cup \left\{ \frac{\sin \frac{m}{n}\pi}{(\frac{m}{n}+2k)\pi}
 : n \ge 2\,,\;1 \le m \le n-1\,,\; k \in \Z \right\}\qquad (q=r)\,.
$$
In \cite{BaWo}, an orthonormal basis of $\ell^2(\DL^0)$ consisting of eigenvectors 
of $P$ is computed. One can of course adapt Corollary~\ref{cor:finitecase} 
in order to transfer the latter into
an orthonormal basis of the subspace $\MC$ of $L^2(\DL^1)$ consisting of 
eigenvectors of $A$. 
\end{exa}

\end{document}